\pdfoutput=1
\documentclass[12pt]{amsart}
\usepackage{amsmath,amsfonts,amssymb,amsthm,amsxtra,pstricks}
\usepackage{hyperref}
\usepackage[english]{babel}         %%wegen Datumsangabe!

\usepackage{epsfig,epic,eepic,graphics,rotating,graphicx,bbm,texdraw,color}
\usepackage[all,cmtip]{xy}
\usepackage{slashbox}

\newcommand{\M}{{\mathbb M}}

\newcommand{\A}{{\mathbb A}}

\newcommand{\C}{{\mathbb C}}

\renewcommand{\P}{{\mathbb P}}

\newcommand{\F}{{\mathbb F}}
\newcommand{\bF}{{\boldsymbol{\mathbb F}}}

\newcommand{\bm}{{\mathbf M}}
\newcommand{\bt}{{\mathbf T}}
\newcommand{\bx}{{\mathbf X}}

\newcommand{\be}{{\mathbf E}}
\newcommand{\bbf}{{\mathbf F}}

\newcommand{\ke}{{\mathcal E}}
\newcommand{\bke}{{\boldsymbol{\mathcal E}}}
\newcommand{\kf}{{\mathcal F}}
\newcommand{\bkf}{{\boldsymbol{\mathcal F}}}

\newcommand{\ki}{{\mathcal I}}

\newcommand{\ko}{{\mathcal O}}

\newcommand{\kq}{{\mathcal Q}}

\newcommand{\kt}{{\mathcal T}}
\newcommand{\ku}{{\mathcal U}}

\DeclareMathOperator{\id}{id}

\DeclareMathOperator{\SL}{SL}

\DeclareMathOperator{\Hom}{Hom}
\DeclareMathOperator{\Aut}{Aut}

\DeclareMathOperator{\Span}{Span}

\newcommand{\ttop}{T_{top}}

\newcommand{\lra}{\longrightarrow}

\newcommand{\xra}{\xrightarrow}

\newcommand{\intoo}[1]{\:
\xymatrix@1{\ar@{^(->}[r]^{#1}&}\:
}

\begin{document}
\newtheorem{sub}{}[section]
\newtheorem{subsub}{}[sub]
\newtheorem{subsubsub}{}[sub]

\title[Trees of vector bundles]{A Kirwan blow-up and trees of vector 
bundles}

\author[Trautmann]{G.~Trautmann}
\address{Universit\"at Kaiserslautern, Fachbereich Mathematik,
\newline  Erwin-Schr\"odinger-Stra{\ss}e
\newline D-67663 Kaiserslautern}
\email{trm@mathematik.uni-kl.de}


\begin{abstract}
In the paper \cite{MTT} a conceptuel description of compactifications
of moduli spaces of stable vector bundles on surfaces has been given,
whose boundaries consist of vector bundles on trees of sufaces.
In this article a typical basic case for the projective plane
is described explicitly including the constrution of a relevant Kirwan
blow up.

%$$\bbf, \abF, \F, \bF  \bkf$$

\end{abstract}

\vspace{2cm}
\maketitle

\begin{section}{Introduction}  To some extent, the replacement 
of limit sheaves in a compactification of a space of vector bundles by 
vector bundles on trees of surfaces is very natural, being in analogy 
to bubbling phenomena in Geometric Analysis and Yang-Mills theory in the
work of Taubes, Uhlenbeck and Feehan. There the degeneration of connections
and fields is described by a process where data are preserved by shifting
them partially to a system of attached 4-spheres. In the analogous
situation of algebraic moduli spaces of vector bundles the attached 
4-spheres can be replaced by projective planes $\P_2$ hanged in at
exceptional lines after blowing up points in a given surface. Then a limit 
sheaf
can be transformed eventually to a vector bundle on the new reducible surface
or on repeatedly constructed trees of surfaces. In \cite{MTT} the trees 
of surfaces and vector bundles have been defined so that these objects
can be the points of a compactification of the moduli spaces of rank-2
vector bundles on a given algebraic surface and are minimal for that purpose.
The original basic example of such a compactification is the moduli space
$M(2;0,2)$ of stable rank-2 vector bundles with Chern classes $c_1=0, c_2=2$
on $\P_2$ which has partially been treated in \cite{MTT}. In this paper 
an explicit construction of the Kirwan blow up of a relevant parameter space is 
given together with the construction of a universal family. In section 2
we recall shortly the definitions and the main theorem of \cite{MTT} 
and in section 3 the typical limit trees are explicitly constructed.

Notation: All varieties in this article shall be defined over an algebraically 
closed field $k$ of characteristic zero. $P(V)$ denotes the projective space
of lines in the $k$-vector space $V$, whereas $\P_n=P(k^{n+1}).$
The points of $P(V)$ are written as [v].

\end{section}

\begin{section}{Trees of surfaces and bundles}

\begin{sub}\label{trs}{\bf Trees.}\rm

A {\bf tree} $T$ in this article is a finite graph, oriented by a partial
order $\le $ and satisfying:
\begin{itemize}
\item there is a unique minimal vertex $\alpha\in T$, the root of $T$;
\item for any $a\in T,\; a\neq \alpha$, there is a unique maximal vertex
$b<a$, the predecessor of $a$, denoted by $a^-$;
\item By $a^+:=\{b\in T\ |\ b^- =a\}$ we denote the set of direct
successors of $a\in T$. We let $\ttop$ denote the vertices of $T$
without successor.
\end{itemize}

A tree of surfaces over a given smooth projective surface $S$, 
modelled by a tree $T$, is a union
$$S_T=S_\alpha\cup\bigcup_a S_a$$
where 
\begin{itemize}
\item $S_\alpha$ is a blow-up of $S$ in finitely many points
\item for $a\in\ttop$, $S_a$ is a projective plane $P_a=\P_2$
\item if $\alpha\not=a\not\in\ttop$, $S_a$ is a blown up projective 
plane $P_a=\P_2$
in finitely many simple points not on a line $l_a\subset P_a$
\item if $a\not=\alpha$, $S_a\cap S_{a^-}=l_a$ and $l_a$ is an exceptional 
line in $S_{a^-}$
\end{itemize}

Such trees can be construted by consecutive blow-ups of simple points,
hanging in a $\P_2(k)$ in each exceptional line of the previous surface and 
then blowing up points in the new $\P_2$, the whole starting with the 
given surface $S$. 

By the construction of $S_T$, all or a part of its components
can be contracted. In particular, there is the morphism
$$S_T\xra{\sigma} S$$
which contracts all the components except $S_\alpha$ to the points of 
the blown up finite set of $S_\alpha$.

Note that:\\
1) There are no intersections of the components other than the lines $l_a$.\\
2) If $T=\{\alpha\}$ is trivial, then $S_T=S.$\\
3) After contracting the lines $l_a$ topologically (when defined over $\C$), 
one obtains bubbles of attached $4$-spheres.

\end{sub}

\begin{sub}\label{trvb}{\bf Treelike vector bundles.}\rm

A {\it weighted tree} is a pair $(T, c)$  of a tree $T$ with a map $c$
which assigns to each vertex $a\in T$ an integer $n_a\geq 0$, called
the {\it weight} or {\it charge} of the vertex, subject to 
$$\#a^+\geq 2\ \text{ if }\ n_a=0 \text{ and } a\ne\alpha.$$
The total weight or total charge of a weighted tree is the sum
$\Sigma_{a\in T}\; n_a=n$ of all the weights.
We denote by $\bt_n$ the set of all trees which admit
a weighting of total charge $n$. It is obviously finite.

In the following we consider only pairs $(S_T,E_T)$, called  
{\bf $\bt_n$-bundles} or simply\\ {\bf tree bundles}, where $T\in\bt_n$,
$S_T$ is a tree of surfaces
and $E_T$ is a rank-2 vector bundle on $S_T$, such that
$c_1(E_T|S_a)=0$, $c_2(E_T|S_a)=n_a$ for all weights $n_a$, and such that the 
bundles $E_a=E_T|S_a$ are
{\bf ``admissible''}, replacing a lacking stability condition, see \cite{MTT}.

In case $S_T=S$ this includes that the bundle $E$ on $S$ belongs to 
$M_{S,h}^b(2;0,n)$, the
quasi-projective Gieseker--Maruyama moduli scheme of $\chi$-stable 
rank-2 vector bundles on $S$ with respect to a polarization $h$
and of Chern classes $c_1=0, c_2=n$.
The bundles in the special case of this article will all be  admissible.

In particular, an indecomposable bundle $E_a$ on $P_a=\P_2$ will 
be admissible if $c_1=0, c_2=1$. 
Such a bundle is not semistable on $P_a$. It is represented 
in homogeneous coordinates by exact sequences 
$$0\to\ko_{P_a}(-2)\xra{(z_0^2, z_1, z_2)}\ko_{P_a}\oplus 2\ko_{P_a}(-1)\to E_a\to 0.$$  

We call the so defined tree bundles also $\bt_n$-bundles.
There is a natural notion of isomorphism for the pairs $(S_T,E_T).$
They consist of isomorphisms of the surfaces with the base surface fixed,
and of isomorphisms of the lifted bundles.
\end{sub}

\begin{sub}\label{ftrvb}{\bf Families of tree bundles.}\rm

A  $\bt_n$-family of tree bundles is a triple $(\be/\bx/Y),$
where $\bx$ is flat family of $\bt_n$-surfaces $X_y,\; y\in Y,$ 
and $\be$ is a rank-2 vector bundle on $\bx$ such that each $E_y=\be|X_y$
is a $\bt_n$-bundle.  

One can then consider the moduli stack $\M_n$ defined by
$$\M_n(Y):= \text{ set of families } (\be/\bx/Y)$$ 
such that any bundle $E_y=\be|X_y$ is a 1-parameter limit of bundles
in $M^b_{S,h}(2;0,n)$. Let
$$\bm_n(Y)=\M_n(Y)/\sim.$$
be the associated functor. The following theorem is stated in \cite{MTT}.

{\bf Theorem:} {\it There is a separated algebraic space $M_n(S)$ of 
finite type over $k$ corepresenting the functor $\bm_n$.}

However the following questions are still open:
\begin{itemize}
\item Is $M_n(S)$ complete?
\item When is $M_n(S)$ a (projective) scheme? 
\item Is $M_n(\P_2)$ a projective compactification of $M_{\P_2}(2;0,n)$?
\item Classification of limit tree bundles for $M_{\P_2}(2;0,n)$ for $n\ge 3$?
\item What about higher rank bundles on $\P_2$?
\item Limit treelike bundles for instanton bundles on $\P_3$?
\end{itemize}

\end{sub}
\end{section}

\begin{section}{Limit trees for $M^b(0,2)$}

Let $M(2;0,2)$ be the moduli space of semistable sheaves 
on $\P_2$ with Chern classes $c_1=0, c_2=2$ and rank $2$ and let
$M^b(0,2)$ be its open part of (stable) bundles.
It is well known that
$M(2;0,2)$ is isomorphic to the $\P_5$ of conics in the dual plane, the
isomorphism being given by $[\kf]\leftrightarrow C(\kf),$ where $[\kf]$
is the isomorphism class of $\kf$ and $C(\kf)$
is the conic of jumping lines of $[\kf]$ in the dual plane. 

It is also well known that any sheaf $\kf$ from $M(2;0,2)$ has
two Beilinson resolutions on $P=\P_2=P(V)$
\begin{equation}\label{beilin}
0\to 2\, \Omega^2_P(2)
\overset{A}{\lra} 2\, \Omega_P^1(1)\to\kf\to 0
\end{equation}
$$0\to 2\, \ko_P(-2)
\overset{B}{\lra} 4\, \ko_P(-1)\to\kf\to 0,
$$
where the matrices $A$ (of vectors in $V$) and $B$ (of vectors in $V^*$)
are related by the exact sequence
$$
0\to k^2\overset{A}{\lra}k^2\otimes V \overset{B}{\lra}k^4\to 0.
$$

The conic  $C(\kf)$ in the dual plane has the equation $\det(A)$.

$\kf$ is locally free if and only if $C(\kf)$ is smooth or if and only if 
$\kf$ is stable. If $C(\kf)$ decomposes into a pair of lines,
then $A$ is equivalent to a matrix of the form
$\left(\begin{smallmatrix}x& 0\\z & y\end{smallmatrix}\right),$
and then $\kf$ is an extension
\[
0\to \ki_{[x]}\to \kf\to\ki_{[y]}\to 0,
\]
whose extension class is represented by the entry $z$. 

{\bf Notice} here that the sheaf is still locally free at the point [y]
if the extension class is non-zero, i.e. $z\not\in\Span(x,y).$
In any case $\kf$ is $S$-equivalent to the direct sum
$\ki_{[x]}\oplus\ki_{[y]}$ .

\begin{sub}\label{typ1}{\bf Type 1 degeneration}\rm

In the following let $e_0,\; e_1,\; e_2$ be basis of $V$ and denote
by $x_0,\; x_1,\; x_2$ its dual basis. For the first example, consider the 
1-parameter deformation 
$\left(\begin{smallmatrix}e_0 & tae_1\\tbe_2 & e_0 \end{smallmatrix}\right)$
with second Beilinson resolution
$$0\to 2\, \ko_C\boxtimes\ko_P(-2)
\overset{B(t)}{\lra} 4\,\ko_C\boxtimes\ko_P(-1)\to\F\to 0,$$
$$B(t)=\left(\begin{smallmatrix}x_1 & x_2 & tax_0 & 0\\
0 & tbx_0 & x_1 & x_2\end{smallmatrix}\right)$$

with parameters $a,b$, where $C=\A^1(k)$. For $t=0$ the sheaf $\F_0$
is singular at $p=[e_0]$, $\F_0=\ki_p\oplus\ki_p.$
The blowing-up $\sigma :Z\to C\times P$ at $(0,p)$
is the subvariety of $C\times P\times \P_2$
given by the equations
$$tx_0u_1-x_1u_0 = 0,\quad tx_0u_2-x_2u_0 = 0,\quad x_1u_2-x_2u_1 = 0,$$
where the $u_\nu$ are the coordinates of the third factor $\P_2.$
We consider the following divisors on $Z$:
\begin{itemize}
\item $\tilde{P}$, the proper transform of $\{0\}\times P$, 
isomorphic to the blow-up of $P$ at $p$; 
\item $D$, the exceptional divisor of $\sigma$; 
\item $H$, the lift of $C\times h$, where $h$ is a general line in $P$;
\item $F$, the divisor defined by $\ko_Z(F)=pr_3^*\ko_{\P_2}(1)$,
\end{itemize}
as shown in the figure 

\hspace*{2cm}\scalebox{0.5}{\begin{picture}(0,0)%
\includegraphics{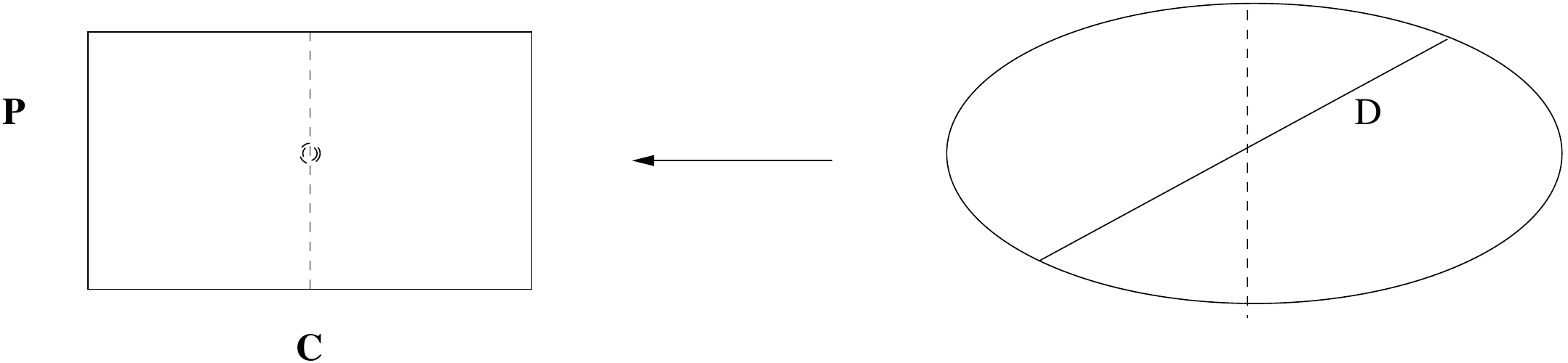}%
\end{picture}%
\setlength{\unitlength}{4144sp}%
\begingroup\makeatletter\ifx\SetFigFont\undefined%
\gdef\SetFigFont#1#2#3#4#5{%
  \reset@font\fontsize{#1}{#2pt}%
  \fontfamily{#3}\fontseries{#4}\fontshape{#5}%
  \selectfont}%
\fi\endgroup%
\begin{picture}(9858,2272)(1561,-2911)
\put(9136,-916){\makebox(0,0)[lb]{\smash{{\SetFigFont{14}{16.8}{\rmdefault}{\mddefault}{\updefault}{\color[rgb]{0,0,0}$\tilde{P}$}%
}}}}
\put(11251,-2221){\makebox(0,0)[lb]{\smash{{\SetFigFont{14}{16.8}{\rmdefault}{\mddefault}{\updefault}{\color[rgb]{0,0,0}Z}%
}}}}
\put(5986,-1591){\makebox(0,0)[lb]{\smash{{\SetFigFont{14}{16.8}{\rmdefault}{\mddefault}{\updefault}{\color[rgb]{0,0,0}$\sigma$}%
}}}}
\end{picture}%
}

Then $D\sim H-F$, and we let $x_\nu$ resp. $u_\nu$ denote the sections
of $\ko_X(H)$ resp. $\ko_Z(F)$ lifting the above coodinates. Using the
equations of $Z$, we see that the canonical
section $s$ of $\ko_Z(D)$ is a divisor of the sections $x_\nu$, such that
$tx_0=su_0$, $x_1=su_1$, $x_2=su_2$ and gives rise to the diagram
{\scriptsize
\[
\xymatrix{ 
 &  & & 2\,\ko_D(-1)\ar @{>->}[d] & \\
0\ar[r]& 2\,\ko_Z(-2H)\ar[r]^{\sigma^*B(t)}\ar @{>->}[d]^s & 4\,\ko_Z(-H)\ar[r]
\ar@{=}[d] & \sigma^*\F\ar[r]\ar @{>->>}[d] & 0\\
0\ar[r] & 2\,\ko_Z(-H-F)\ar[r]^(0.6){B_Z}\ar @{>->>}[d] & 4\,\ko_Z(-H)\ar[r]& 
\bbf\ar[r]  & 0\\ & 2\,\ko_D(-1) &  &   & }
\]}

with $B_Z=\left(\begin{smallmatrix}u_1 & u_2 & au_0 & 0\\
0 & bu_0 & u_1 & u_2\end{smallmatrix}\right).$ 
Thus $B_Z$ represents
a locally free sheaf $\bbf$ on $Z$, but its first Chern class
has been modified by blowing up and removing the torsion. 
To correct this, consider the twisted bundle $\be:=\bbf(D)$.
Then 
$\be|{\tilde{P}}\simeq 2\,\ko_{\tilde{P}},$ and the restriction
$\be|D$ belongs to $M^b_{D}(2;0,2)$, $D\simeq\P_2.$ Moreover, $Z$ is flat over 
$C$ and $\be$ is a flat family of vector bundles over $C$ with the limit tree
bundle $\be|Z_0$ on the fibre $Z_0=\tilde{P}\cup D$ over $0\in C.$
This can be symbolized by 

\hspace*{5cm}\scalebox{0.5}{\begin{picture}(0,0)%
\includegraphics{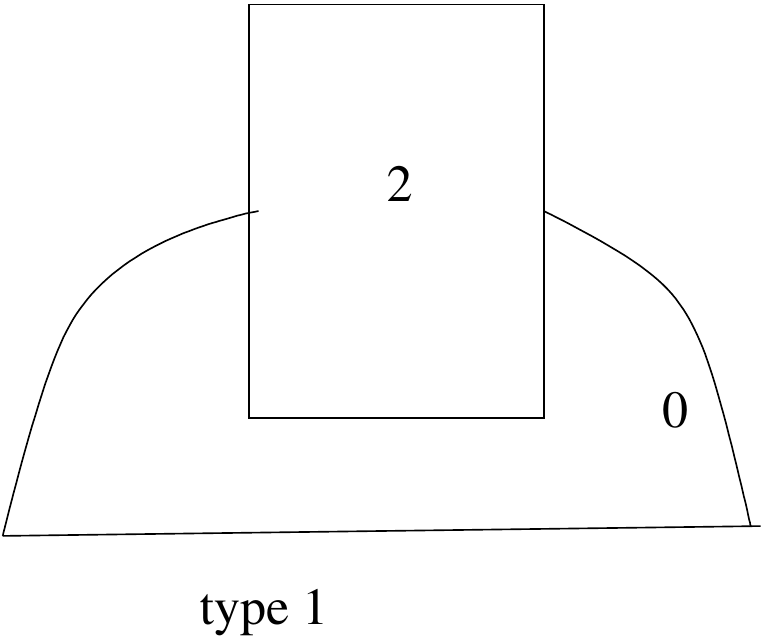}%
\end{picture}%
\setlength{\unitlength}{4144sp}%
\begingroup\makeatletter\ifx\SetFigFont\undefined%
\gdef\SetFigFont#1#2#3#4#5{%
  \reset@font\fontsize{#1}{#2pt}%
  \fontfamily{#3}\fontseries{#4}\fontshape{#5}%
  \selectfont}%
\fi\endgroup%
\begin{picture}(3489,2929)(3859,-3743)
\end{picture}%
}

the numbers indicating the second Chern classes of the bundles on 
the components.
The isomorphism class of this limit depends on the chosen parameters
$a,b$ which determine a normal direction to the Veronese surface
in $\P_5$. This leads to blowing it up and to the Kirwan blow-up
of the parameter space, see section 4.

\end{sub}

\begin{sub}\label{typ2}{\bf Type 2 degeneration}\rm

Let now a family on $C\times P$ be given by 
$\left(\begin{smallmatrix}e_0 & -te_1\\-te_1 & e_2 \end{smallmatrix}\right),$
defining a deformation of the sheaf of 
$\left(\begin{smallmatrix}e_0 & 0\\0& e_2 \end{smallmatrix}\right).$
Similarly to the previous case, the deforming sheaf $\F$ is the cokernel
of the matrix
$$B(t)=\left(\begin{matrix}x_2 & x_1 & tx_0 & 0\\
0 & tx_2 & x_1 & x_0\end{matrix}\right)$$

Blowing up $C\times P$ in the two singular points $(0,p_0)$ and $(0,p_2)$,
$p_\nu=[e_\nu]$, leads to the figure

\hspace*{3cm}\scalebox{0.5}{\begin{picture}(0,0)%
\includegraphics{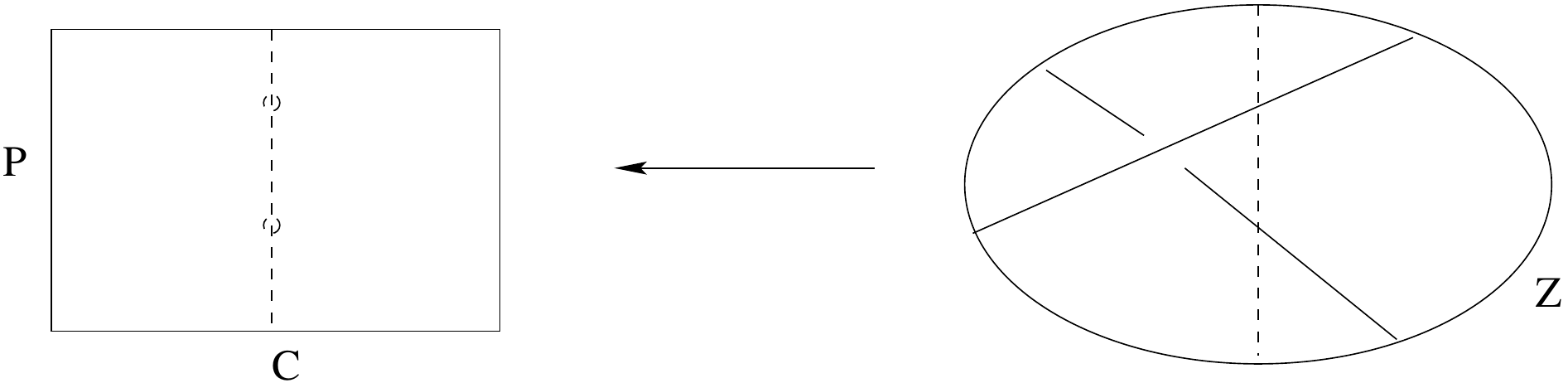}%
\end{picture}%
\setlength{\unitlength}{4144sp}%
\begingroup\makeatletter\ifx\SetFigFont\undefined%
\gdef\SetFigFont#1#2#3#4#5{%
  \reset@font\fontsize{#1}{#2pt}%
  \fontfamily{#3}\fontseries{#4}\fontshape{#5}%
  \selectfont}%
\fi\endgroup%
\begin{picture}(8643,2097)(1336,-2416)
\put(8731,-916){\makebox(0,0)[lb]{\smash{{\SetFigFont{14}{16.8}{\rmdefault}{\mddefault}{\updefault}{\color[rgb]{0,0,0}$D_0$}%
}}}}
\put(8641,-1771){\makebox(0,0)[lb]{\smash{{\SetFigFont{14}{16.8}{\rmdefault}{\mddefault}{\updefault}{\color[rgb]{0,0,0}$D_2$}%
}}}}
\put(5221,-1141){\makebox(0,0)[lb]{\smash{{\SetFigFont{14}{16.8}{\rmdefault}{\mddefault}{\updefault}{\color[rgb]{0,0,0}$\sigma$}%
}}}}
\put(7921,-601){\makebox(0,0)[lb]{\smash{{\SetFigFont{14}{16.8}{\rmdefault}{\mddefault}{\updefault}{\color[rgb]{0,0,0}$\tilde{P}$}%
}}}}
\end{picture}%
}\\

The blown up variety $Z$ has the standard embedding into 
$(C\times P)\times\P_2\times\P_2$ with divisors 
\begin{itemize}
\item $H$, the pull back of the divisor $C\times h$ in $C\times P$\\
\item $\tilde{P}$, the blow-up of $\{0\}\times P$ in the two points\\
\item $D_0$, $D_2$, the two exceptional divisors\\
\item $F_0$, $F_2$, whose invertible sheaves are the pull backs of 
$\ko_{\P_2}(1)$ from the third and fourth factor.
\end{itemize}
Letting denote $x_\nu, u_\nu, v_\nu$, and $s_0, s_2$ the basic sections of
the sheaves of $H, F_0, F_2, D_0, D_2$, we have the equations (as homomorphisms
between invertible sheaves) 
$tx_0=s_0u_0,\quad x_1=s_0u_1,\quad x_2=s_0u_2,$ and 
$x_0=s_2v_0,\quad x_1=s_2v_1,\quad tx_2=s_2v_2.$

By that we have the matrix decomposition
$$\left(\begin{matrix} x_2 & x_1 & tx_0 & 0\\
0 &tx_2 & x_1 & x_0  \end{matrix}\right)\left(\begin{matrix}s_0 & 0\\
0 & s_2 \end{matrix}\right) =
\left(\begin{matrix} u_2 & u_1 & u_0 & 0\\
0 & v_2 & v_1 & v_0\end{matrix}\right).$$
Using this, the torsion of $\sigma^*\F$ can be removed as in the diagram of 
the previous section. Then  $\bbf=\sigma^*\F/torsion$ has the resolution
$$0\to\ko_Z(-H-F_0)\oplus\ko_Z(-H-F_2)\xra{B_Z} 4\ko_Z(-H)\to \bbf\to 0,$$
where $B_Z$ is the right hand matrix. The tree components of $\bbf$ are
$\bbf|D_i=\kt_{D_i}(-1),$
$\bbf|\tilde{P}=\ko_{\tilde{P}}(-l_0)\oplus\ko_{\tilde{P}}(-l_2),$
where $l_0, l_2$ are the exceptional lines on $\tilde{P}.$
However, there is no way by twist or elementary transformation to make
the first Chern classes $c_1$ vanish.

But starting with 
$\left(\begin{smallmatrix}e_0 & -t^2e_1\\-t^2e_1 & e_2 \end{smallmatrix}\right),$
we get by the same procedure a sheaf $\bbf$ on $Z$ whose resolution
matrix is 
$$B_Z=\left(\begin{matrix} u_2 & u_1 & tu_0 & 0\\
0 & tv_2 & v_1 & v_0\end{matrix}\right).$$
This resolution implies
that $\bbf$ is reflexive and singular in exactly two points 
$q_0=\{u_1=u_2=t=0\}$ and $q_2=\{v_1=v_0=t=0\}$, and that its restrictions
to the components of $Z_0=\tilde{P}\cup D_0\cup D_2$ are
 
$\bbf|\tilde{P}=\ko_{\tilde{P}}(-l_0)\oplus\ko_{\tilde{P}}(-l_2)$ and
$\bbf|D_i=\ko_{D_i}\oplus\ki_{q_i,D_i}(1).$

Hence there is an elementary transform on Z,
$$0\to\bbf^\prime\to\bbf\to\ko_{D_0}\oplus\ko_{D_2}\to 0.$$
The resolution of $\bbf^\prime$ can be computed as follows. There is a
decomposition $tu_0=s_0\tilde{u_0}$ because $tu_0$ vanishes on the divisor
$D_0$. Similarly we have 
$tv_2=s_2\tilde{v_2}, \; u_1=s_2\bar{u_1},\; v_1=s_0\bar{v_1},$
and from this the matrix decomposition 
$$\left(\begin{matrix} u_2 & \bar{u_1} & \tilde{u_0} & 0\\
0 & \tilde{v_2} & \bar{v_1} & v_0 \end{matrix}\right)
\left(\begin{smallmatrix} 1 & & & \\
 & s_2 &  &  \\
 &  & s_0 &  \\
 &  &  & 1\end{smallmatrix}\right) =
\left(\begin{matrix} u_2 & u_1 & tu_0 & 0\\
0 & tv_2 & v_1 & v_0\end{matrix}\right)$$

It follows by diagram chasing that the left hand matrix gives the resolution
$$0\to\ke_1\to\ke_0\to\bbf^\prime\to 0,$$ 
where $\ke_1=\ko_Z(-H-F_0)\oplus\ko_Z(-H-F_2)$ and\\ 
$\ke_0=\ko_Z(-H)\oplus\ko_Z(-H-S_2)\oplus\ko_Z(-H-S_0)\oplus\ko_Z(-H)$

This resolution shows that $\bbf^\prime$ is locally free on $Z$.
In order to determine its restrictions to the components, one should
use the identities $u_0^2=x_0\tilde{u_0}$, $v^2_2=x_2\tilde{v_2}$,
which follow from the previous identities. Using these, one can determine
the restrictions of the twisted bundle $\be:=\bbf^\prime(D_0+D_2):$

$\be|\tilde{P}=2\ko_{\tilde{P}}$ and
$\be|D_i$ is a bundle on $D_i\cong\P_2$ with Chern classes $c_1=0, c_2=1$ 
(see the description of bundles with these Chern classes in \ref{trvb}.)

Since the elementary transform and the twisting do not affect the bundle
on the part of $Z$ over $C\setminus\{0\}$, the sheaf is a limit tree bundle
on the fibre $Z_0=\tilde{P}\cup D_0\cup D_2.$

\end{sub}

\begin{sub}\label{typ3}{\bf Type 3 degeneration}\rm

Let $\tilde{\P}_5$ be the blow-up of $P(S^2V)=\P_5$ of the Veronese
surface in $\P_5$, let $\Sigma_2\subset\tilde{\P}_5$ be the exceptional
divisor and $\Sigma_1\subset\tilde{\P}_5$ the proper transform
of the divisor of degenerate conics, see also \ref{veronese}.

By the above, type 1 limit tree bundles belong to $\Sigma_2\setminus\Sigma_1$
and type 2 limit tree bundles belong to $\Sigma_1\setminus\Sigma_2.$
There is a third type of limit tree bundle belonging to 
$\Sigma_2\cap\Sigma_1$ with symbolic tree

\hspace*{5cm}\scalebox{0.5}{\begin{picture}(0,0)%
\includegraphics{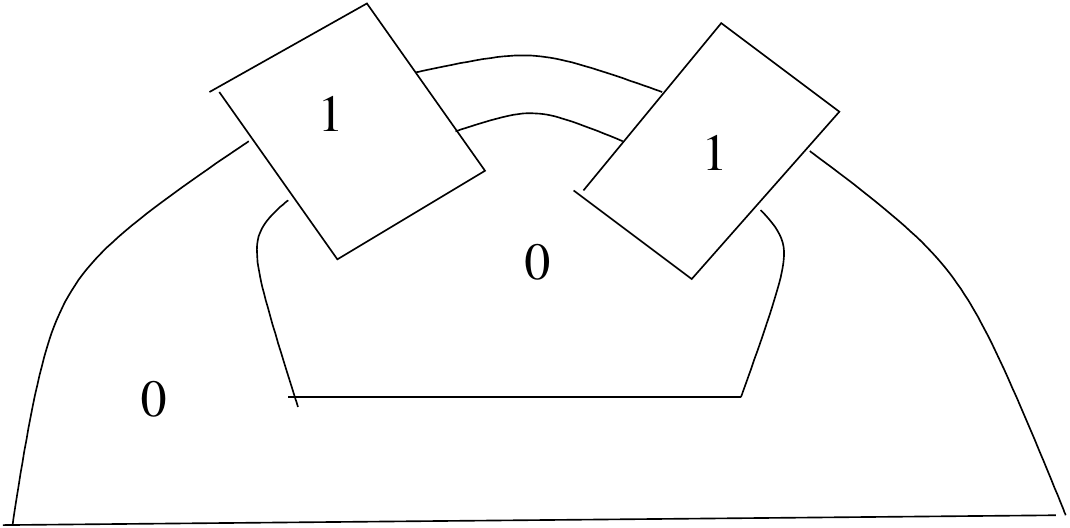}%
\end{picture}%
\setlength{\unitlength}{4144sp}%
\begingroup\makeatletter\ifx\SetFigFont\undefined%
\gdef\SetFigFont#1#2#3#4#5{%
  \reset@font\fontsize{#1}{#2pt}%
  \fontfamily{#3}\fontseries{#4}\fontshape{#5}%
  \selectfont}%
\fi\endgroup%
\begin{picture}(4884,2409)(3499,-2593)
\end{picture}%
}

Examples can be obtained as limits of families of type 
$\left(\begin{smallmatrix}e_0 & -t^3e_1\\-t^3e_1 & e_0+te_2\end{smallmatrix}
\right)$
and two consecutive blow-ups. In this case the family $\F$ on $C\times P$ 
is given as the cokernel in
$$0\to 2\, \ko_C\boxtimes\ko_P(-2)
\overset{B(t)}{\lra} 4\,\ko_C\boxtimes\ko_P(-1)\to\F\to 0,$$
$$B(t)=\left(\begin{smallmatrix}x_2 & x_1 & t^3x_0 & 0\\
  0 & t^2x_2 & x_1 & tx_0-x_2\end{smallmatrix}\right).$$
This sheaf $\F$ is singular in $(0,p)$, $p=[e_0]$. Let then
$$\sigma : Z\to  C\times P$$
be the blow-up as in \ref{typ1}, described as subvariety of 
$Z\subset C\times P\times \P_2$ with divisors 
$\tilde{P},\; H,\; D,\; F,\; D\sim H-F$. Let $s$ be the standard section
of  $\ko_Z(D)$ for the exceptional divisor, and let $x_\nu$ respectively $y_\nu$
be the basic sections of $\ko_Z(H)$ respectively $\ko_Z(F)$ with equations
$tx_0=sy_0,\; x_1=sy_1,\; x_2=sy_2.$ It follows as in \ref{typ1} that
the sheaf $\bbf=\sigma^*\F/torsion$ has the resolution
$$0\to 2\ko_Z(-H-F)\xra{B_Z} 4\ko_Z(-H)\to \bbf\to 0,$$
$$B_Z=\left(\begin{smallmatrix}y_2 & y_1 & t^2y_0 & 0\\
0 & t^2y_2 & y_1 & y_0-y_2\end{smallmatrix}\right).$$ 
This sheaf and its syzygy is of the same type as in \ref{typ2}.
It is reflexive and singular exactly in the points 
$p_0,\; p_2\in D\smallsetminus\tilde{P}$,\quad
$p_0=\{t=y_1=y_2=0\}$ and $p_2=\{t=y_1=y_0-y_2=0\}.$ 
Againe one can verify that the sheaf $\bbf^\prime:=\bbf(D)$ has the restictions
$$\bbf^\prime|\tilde{P}=2\ko_{\tilde{P}}\quad\text{and}\quad \bbf^\prime|D=\ki_{p_0,D}\oplus\ki_{p_2,D}.$$
on the components of $Z_0=\tilde{P}\cup D.$

In order to construct a locally free limit tree bundle we blow-up $Z$
in the two points $p_0,\; p_2$ to get
$$\tau : W\to Z$$
with exceptional divisors $S_0$ and $S_2$, the proper transform $\tilde{D}$
of $D$, the lifted divisors $\tilde{P}$ and $F$, and the two divisors 
$F_0$ and $F_2$ coming from the embedding.

As in \ref{typ2} one concludes that the sheaf 
$\bbf^{\prime\prime}=\tau^*\bbf^\prime/torsion$ is reflexive and the cokernel
of a matrix

$$B_Z=\left(\begin{matrix} u_2 & u_1 & tu_0 & 0\\
0 & tv_2 & v_1 & v_0\end{matrix}\right),$$

and such that $\bbf^{\prime\prime}$ restricts as
$$\bbf^{\prime\prime}|\tilde{P}=2\ko_{\tilde{P}},\quad \bbf^{\prime\prime}|\tilde{D}=\ko_{\tilde{D}}(-l_0)\oplus\ko_{\tilde{D}}(-l_2),\quad \bbf^{\prime\prime}|S_i=\ko_{S_i}\oplus\ki_{q_0,S_i}(1),$$
where $q_i\in S_i\smallsetminus\tilde{D}.$

Finally, as in \ref{typ2}, there is an elementary transform
$$0\to\be^\prime\to\bbf^{\prime\prime}\to\ko_{S_0}\oplus\ko_{S_0},$$
such that $\be^\prime$ is locally free on $W$ and such that
$\be := \be^\prime(S_0+S_2)$ has the desired restrictions
$$\be|\tilde{P}=2\ko_{\tilde{P}},\quad \be|\tilde{D}=2\ko_{\tilde{D}},\quad$$
and such that $\be|S_i$ do have the Chern classe $c_1=0,\; c_2=1.$
So $\be$ is a limit tree bundle on the tree of surfaces 
$W_0=\tilde{P}\cup\tilde{D}\cup S_0\cup S_2.$
\end{sub}
\end{section}

\begin{section}{Kirwan blow-up I}
The $2\times 2$-matrices with entries in $V$ in (\ref{beilin}) parametrize the 
sheaves in $M(2;0,2)$ and at the same time the conics of their jumping lines
in the dual plane $P(V^*)$ by their determinants in $S^2V.$ Since the 
isomorphisms of the left hand term in (\ref{beilin}) are not essential, only
the subspaces $[A]$ spanned by the rows of $A$ matter, so that
the Grassmannian $G_2(k^2\otimes V)$ is a parameter space of $M(2;0,2).$
The Pl\"ucker embedding 
$$ p: G_2(k^2\otimes V)\subset P(\wedge^2(k^2\otimes V))=
P(\wedge^2k^2\otimes S²V\oplus S^2k^2\otimes\wedge^2V)$$
can be expressed in terms of the entries, using the standard basis of $k^2$, by
$$[\left(\begin{matrix} x & x^\prime\\ y & y^\prime \end{matrix}\right)]
\quad \overset{p}\mapsto\quad [xy^\prime-x^\prime y;\; x\wedge y,\; x\wedge y^\prime + x^\prime\wedge y,\;  x^\prime\wedge y^\prime].$$ 
One should note here that there is the relation 
\begin{equation}\label{wedgerel}
\left(\begin{matrix} x & x^\prime\\ y & y^\prime \end{matrix}\right)\wedge
\left(\begin{matrix} x\wedge y & x\wedge y^\prime + x^\prime\wedge y & x^\prime\wedge y^\prime & 0\\
0 & x\wedge y & x\wedge y^\prime + x^\prime\wedge y & x^\prime\wedge y^\prime
\end{matrix}\right)=0.
\end{equation}
There is an action of $\SL_2(k)$ on both sides of the Pl\"ucker embedding,
induced by the natural action on $k^2$ and written as 
$$ [A]g=[Ag]\quad\text{and}\quad [q;\; \Phi]g=[q;\; \Phi S^2g],$$
explicitly with
$$Ag=\left(\begin{matrix}x & x^\prime\\ y & y^\prime\end{matrix}\right)
\left(\begin{matrix} \alpha & \beta \\ \gamma & \delta \end{matrix}\right)
\quad\text{and}\quad 
\Phi S^2g=(\xi,\; \omega,\; \eta)\left(\begin{matrix} \alpha^2 & 2\alpha\beta & \beta^2 \\ 
2\alpha\gamma & \alpha\delta + \beta\gamma & 2\beta\delta \\
\gamma^2 & 2\gamma\delta & \delta^2\end{matrix}\right),$$
such that the Pl\"ucker embedding is equivariant. An element [A] in the 
Grassmannian is semistable if and only if $\det(A)\ne 0,$ and it is stable
if and only if $\det(A)$ is the equation of a non-degenerate quadric in the 
dual plane $P(V^*)$. Moreover, the morphism $[A] \to [\det(A)]$ ,
$$G_2(k^2\otimes V)^{ss}\lra PS^2V\cong\P_5\cong\bar{M}(2; 0,2,0)$$
is a good GIT quotient, see \cite{NT}.

For the construction of a compactification of $M^b(0,2)$ by tree bundles we
need to replace the Grassmannian by a parameter space with only stable
points in order to avoid unnatural identifications in the boundary.
This is done by the method of F. Kirwan, \cite{Ki} in two consecutive
blow-ups.

{\bf The first blow-up:}
In the following we use the abbreviations $X=G_2(k^2\otimes V)$ and 
$G=\SL_2(k)$. The group $G$ has the fixed points 
$[\left(\begin{smallmatrix} x & 0\\
0 & x \end{smallmatrix}\right)].$ According to \cite{Ki}, let then $Z_G$ 
denote the subset
$$ Z_G = \{[A]\in X|\; \text{the affine fibre of}\;p(A)\;\text{fixed by}\; G\}.$$ 
It follows that 
$$ Z_G = \{[\left(\begin{smallmatrix} x & 0\\0 & x \end{smallmatrix}\right)]\}
\cong P(V),$$
that it is a closed and smooth subvariety of $X$ and that $Z_G=GZ_G.$ 
The vanishing of the components of $\Phi$ characterizes the points of $Z_G$
and these components define its ideal sheaf $\ki_G$. Let then 
$$\tilde{X}:=Bl_{Z_G}(X)$$
be the blow-up of $X$ along $Z_G$. In this situation 
$$\tilde{X}\subset X\times P(S^2k^2\otimes\wedge^2V)$$ is the closure
of the graph of the map $\Phi: X\smallsetminus Z_G\to P(S^2k^2\otimes\wedge^2V).$ This blow-up can geometrically be described as follows.

\begin{sub}\label{geombl1}{\bf Lemma:}\quad
(a) $\tilde{X}$ is the subvariety of $X\times P(S^2k^2\otimes\wedge^2V)$
of points $([A],\; [\xi,\; \omega,\; \eta])$ satisfying
\begin{enumerate}
\item[(i)] $ (x\wedge y,\; x\wedge y^\prime + x^\prime\wedge y,\; x^\prime\wedge y^\prime)\in k(\xi,\; \omega,\; \eta)$\\
\item[(ii)] $A\wedge\left(\begin{smallmatrix}\xi & \omega & \eta & 0\\
0 & \xi & \omega & \eta \end{smallmatrix}\right)=0$ 
\end{enumerate}
(b) The exceptional divisor $E_G$ in $\tilde{X}$ is the subvariety
of pairs $([A],\; [x\wedge u,\; x\wedge w,\; x\wedge v])$ with 
$A=\left(\begin{smallmatrix} x & 0\\ 0 & x  \end{smallmatrix}\right).$

(c) $\tilde{X}$ is smooth and the projection $\tilde{X}\to X$ is
$G$-equivariant.
\end{sub}

Sketch of proof: Because $\tilde{X}$ is the closure of graph, (i) follows
immediately, and also (ii) by formula (\ref{wedgerel}). Let conversely
$Y\subset X\times P(S^2k^2\otimes\wedge^2V)$ be defined by (i) and (ii).
Then $\tilde{X}\subset Y$ and $\tilde{X}\smallsetminus E_G=Y\smallsetminus E_G.$
One shows now that the fibre $Y_p$ for a point $p\in Z_G$ coincides
with the fibre $\tilde{X}_p=E_{G,p}.$ Such a point has as its first
component 
$A=\left(\begin{smallmatrix} x & 0\\ 0 & x  \end{smallmatrix}\right),$
and (ii) implies that its second component is of the form
$[x\wedge u,\; x\wedge w,\; x\wedge v].$ Consider then the 1-parameter
family $A(t)=\left(\begin{smallmatrix} x & -tv\\ tu & x + tw  \end{smallmatrix}\right).$ For $t\ne 0$, $[A(t)]\in X\smallsetminus Z_G,$ and its lift
to $\tilde{X}$ has the limit $\tilde{A}$ with components 
$[A]$ and $[x\wedge u,\; x\wedge w,\; x\wedge v].$ This proves (a) and also (b)
as a corollary.  For (c) smoothness follows from that of $X$ and $Z_G$,
and the equivariance directly from (a).
%\newpage

It follows from (b) that $E_G$ is the $\P_5$-bundle
\[
\xymatrix{E_G\ar[d]\ar[r]^-{\approx} & P(S^2k^2\otimes Q)\ar[d]\\
Z_G\ar[r]^{\approx} & P(V),
}
\]
where $Q$ is the tautological quotient bundle of P(V). 

{\bf Stability in $\tilde{X}$:}\quad
By definition of $\tilde{X}$ there is the Pl\"ucker embedding
$$\tilde{X}\subset 
P((\wedge^2k^2\otimes S²V\oplus S^2k^2\otimes\wedge^2V)\otimes 
(S^2k^2\otimes\wedge^2V))$$ 
and by this the action on $\tilde{X}$ is induced by the obvious linear 
action of $G$ on the ambient projective space.

\begin{sub}\label{propstab}{\bf Proposition:}\quad 
Let $\tilde{A}=([A],\; [\xi,\; \omega,\; \eta])$ be a point of $\tilde{X}.$  
Then
\begin{enumerate}
\item[(i)]  $\tilde{A}$ is semistable if and only if both of 
$\det(A)=xy^\prime -x^\prime y$ and $\omega^2-4\xi\eta$ are non-zero.
\item[(ii)] If $\tilde{A}\not\in E_G$, then $\tilde{A}$ is stable 
if and only if $\pi(\tilde{A})=[A]$ is stable.
\item[(iii)] If $\tilde{A}\in E_G$, then $\tilde{A}$ is stable 
if and only if $\omega^2-4\xi\eta$ is not a square in $S^2(V/k.x)$
\end{enumerate}
\end{sub}

For the proof, notice first that the quadratic forms 
$\det(A)=xy^\prime -x^\prime y$ and $\omega^2-4\xi\eta$ of the components 
of $\tilde{A}$ are invariant under this action. Then the statements
can be canonically verified by either looking for the points
in the affine cone or by using the Mumford criterion for the action
of 1-parameter subgroups. For the latter, the weights can be computed
via the tensor products in the Pl\"ucker space.

Some elemntary calculations with the explicit description of the group
action show:

\begin{sub}\label{remstab}{\bf Lemma:}\quad
Let $\tilde{A}=([A],\; [\xi,\; \omega,\; \eta])$ be a point of $\tilde{X}.$  
Then
\begin{enumerate}
\item[(i)] $\omega^2-4\xi\eta=0$ if and only if there is a $g\in G$ such
that $[\xi,\; \omega,\; \eta]S^2g=[\xi^\prime, 0, 0].$
\item[(ii)] $\omega^2-4\xi\eta$ is square if and only if there is a $g\in G$ 
such that $[\xi,\; \omega,\; \eta]S^2g=[\xi^\prime,\; \omega^\prime,\; \eta^\prime]$ with $\xi^\prime=0$ or $\eta^\prime=0.$
\item[(iii)] $\omega^2-4\xi\eta$ is a product if and only if there is a 
$g\in G$ 
such that $[\xi,\; \omega,\; \eta]S^2g=[\xi^\prime, 0,\; \eta^\prime].$
\end{enumerate}
\end{sub}

Let now $H_0^{ss}\subset H_1^{ss}\subset X^{ss}$ be the subvarities of
points $[A]$ for which $\det(A)$ is a square respectively a product in $S^2V$.
These are the inverse images in $X^{ss}$ of the double lines respectively
pairs of lines in the space $P(S^2V)$ of conics in $P(V^*).$ 
Let $H_0\subset H_1$ be their closures in $X$. By definition 
$Z_G\subset H_0^{ss}.$ Since the matrices $[A]\in H_0^{ss}$ are of type
$[\left(\begin{smallmatrix} x & 0\\ z & x \end{smallmatrix}\right)]g,\;
g\in G,$ one finds that $H_0^{ss}\smallsetminus Z_G$ consists of all
non-closed orbits whose closures meet $Z_G,$ the orbits of the latter being its 
points. Then 
$$ P(V)\cong Z_G=H_0^{ss}//G\subset X^{ss}//G\cong P(S^2V)$$
is the Veronese embedding. Moreover, all the points 
$H_0^{ss}\smallsetminus Z_G$ become unstable in $\tilde{X}$, see Lemma
\ref{unstab}.

Let denote $\tilde{H}_0\subset \tilde{H}_1$ be the proper transforms of 
$H_0\subset H_1$ in $\tilde{X}.$ Then the following holds.

\begin{sub}\label{unstab}{\bf Lemma:}\quad
(a) $\tilde{H}_0\cap \tilde{X}^{ss}=\emptyset $ and 
$E_G\cap\tilde{H}_0= E_G\smallsetminus E_G^{ss}.$

(b) $\tilde{H}_1\cap \tilde{X}^{s}=\emptyset $ and 
$E_G\cap\tilde{H}_1= E_G\smallsetminus E_G^{s}.$
\end{sub}

Sketch of proof: A point in $H_0^{ss}\smallsetminus Z_G$ is equivalent
to a point $[\left(\begin{smallmatrix} x & 0\\ z & x \end{smallmatrix}\right)]$
and this has the second component $[x\wedge z, 0, 0]$ in $\tilde{X}.$
By Remark \ref{remstab} it is not semistable. Then also 
$\omega^2-4\xi\eta=0$ for the limit points. To show that 
$E_G\smallsetminus E_G^{ss}\subset\tilde{H}_0 $ we may assume that a point
$p$ in $E_G\smallsetminus E_G^{ss}$ has the components 
$[\left(\begin{smallmatrix} x & 0\\ 0 & x \end{smallmatrix}\right)],
[x\wedge z, 0, 0].$ As in the proof of Lemma \ref{geombl1}, the family
defined by 
$[\left(\begin{smallmatrix} x & 0\\ tz & x \end{smallmatrix}\right)]$
shows that $p\in\tilde{H}_0.$ This proves (a). The proof of (b) is
analogous.

By the characterization of semistable points, the 
equivariant morphism $\pi: \tilde{X}\to X$ maps (semi-)stable points
to (semi-)stable and gives rise to a morphism
$\tilde{X}^{ss}//G\to X^{ss}//G\cong P(S^2V)$, which is an isomorphism
over the complement of the Veronese surface $Z_G/G.$
Because $\Sigma_2:=E_G^{ss}//G$ becomes the inverse image of $Z_G/G$
and is a divisor, we obtain the 

\begin{sub}\label{veronese}{\bf Proposition:}
$\widetilde{P(S^2V)}:=\tilde{X}^{ss}//G$  is the blow-up  of 
$P(S^2V)$ along the Veronese surface. 
\end{sub}

\begin{sub}\label{dual}{\bf Related geometry of conics:}\rm\quad
For any point $\tilde{A}$ in $\tilde{X}$ the quadratic form $\omega^2-4\xi\eta$
can be seen as an element of $S^2V^*$ because of $\wedge^2V\cong V^*.$
One can then easily verify that for any non-degenerate 
$A=\left(\begin{smallmatrix} x & x^\prime\\
y & y^\prime \end{smallmatrix}\right)$ (for which
$[x\wedge y,\; x\wedge y^\prime + x^\prime\wedge y,\; x^\prime\wedge y^\prime)]
=[\xi,\; \omega,\; \eta])$),
the quadradic form $\omega^2-4\xi\eta$ is the equation of the dual conic
in $P(V)$ of the conic $\{\det(A)=0\}\subset P(V^*)$ (of jumping lines of the
corresponding vector bundle).
Because $\tilde{X}^{ss}$ is defined by $\det(A)\ne 0$ and 
$\omega^2-4\xi\eta\ne 0$ we can define the universal family of conics
$$Q\subset\tilde{X}^{ss}\times P(V)$$
as the subvariety of pairs $(\tilde{A},[v])$
with  $(\omega^2-4\xi\eta)(v)=0.$
If $\tilde{A}\in E_G^{ss}$, i.e. 
$A=\left(\begin{smallmatrix} x & 0\\ 0 & x \end{smallmatrix}\right)$
then $(\xi,\; \omega,\; \eta)=(x\wedge u,\; x\wedge w,\;  x\wedge v)$
and the fibre $Q_{\tilde{A}}$ is a pair of lines through [x] in $P(V)$
or a double line. 

Secondly, the related quadratic form 
$w^2-4uv\in S^2V/k.x$ without the factor $x$ defines two points or a 
double point on the double line  $\{x^2=0\}$ in $P(V^*).$ 

Recalling that the space of {\bf complete conics} in the plane $P(V^*)$
consists of conics, for which the double lines are enriched by two 
points or a double point,  one finds that $\tilde{X}^{ss}$ parametrizes this 
space and that the quotient $\widetilde{P(S^2V)}:=\tilde{X}^{ss}//G$ is the 
space of complete conics in  $P(V^*).$
Moreover, because the forms $\det(A)$ and $\omega^2-4\xi\eta$
are invariant, the conic bundle Q descents to a conic bundle embedded
in $\widetilde{P(S^2V)}\times P(V)$ and describes the duality for complete
conics.
\end{sub}
\end{section}

\begin{section}{Kirwan Blow-up II}\rm

It is easy to see that there are no semistable points in $\tilde{X}$
with a 2-dimensional stabilizer by checking the types of points. But there
are 1-dimensional such stabilizers. For the Kirwan blow-up it is enough
to consider only connected reductive ones. Again by checking the different
types of points, one finds that the only such stabilizers are 
$R=\{\left(\begin{smallmatrix}\alpha & 0\\ 0 & \alpha^{-1}\end{smallmatrix}\right)\}\cong k^*$ and its conjugates. According to \cite{Ki} we consider 
for the center of the blow-up of $\tilde{X}$ the
subvariety $Z_R$ of points $\tilde{A}$ in $\tilde{X}$ which are fixed
by $R$ and such in addition $R$ acts trivially on the affine fibre of 
$\tilde{A}$ in $\wedge^2(k^2\otimes V)\otimes (S^2k^2\otimes\wedge^2V).$
A direct computation shows that
$$Z_R\quad\text{is the set of points}\quad([\left(\begin{smallmatrix} x & 0\\
0 & y \end{smallmatrix}\right)],\; [0,\omega,0])\quad\text{in\quad}\tilde{X}.$$
Then $GZ_R\subset\tilde{X}^{ss}\smallsetminus\tilde{X}^{s}$ and is of 
dimension 6.

By definition $GZ_R\subset\tilde{H}^{ss}_1$, and $GZ_R$ is the
subset of points in $\tilde{H}^{ss}_1$ with closed orbits. The good quotient
$GZ_R//G=\tilde{H}^{ss}_1//G$ is then the proper transform in 
$\widetilde{P(S^2V)}$ of the divisor 
$\Sigma_1$ of products in $P(S^2V).$

\begin{sub}\label{gzr}{\bf Lemma:}\\
(1) The closure $\overline{GZ}_R$ is the subvariety of points 
$([A], [\xi,\; \omega,\; \eta])$ in $\tilde{X}$ for which 
$\xi,\; \omega,\; \eta$ are pairwise linearly dependent in $\wedge^2V.$

(2) $\overline{GZ}_R\cap\tilde{X}^{ss}=GZ_R. $ 

(3) $\overline{GZ}_R$ is smooth.

(4) $\overline{GZ}_R$ and $E_G$ intersect transversily in dimension 5.
\end{sub}

\begin{proof} Let $Y$ be the closed subvariety of $\tilde{X}$ defined by 
the condition in (1). Then $\overline{GZ}_R\subset Y.$ When 
$y\in Y\cap\tilde{X}^{ss}$, then $y=([A], [a\xi,\; b\xi,\; c\xi])$
with $b^2-4ac\ne 0$ and there is a group element $g$ and some $\lambda$
so that $\lambda (a,b,c)=(0,1,0)S^2g$, because $y$ is supposed to be 
semistable. Then $yg^{-1}=([B], [0,\xi,0])$
and thus an element of $Z_R.$ Now $Y\cap\tilde{X}^{ss}=GZ_R.$  If $y$ is 
unstable, there is a group element
$g$ so that $\lambda (a,b,c)=(1,0,0)S^2g.$ 
Then $yg^{-1}=([B], [\xi,0,0])$ and such points are limits of points in 
$GZ_R$: such matrices $B$ can only be of type $[\left(\begin{smallmatrix}x & 0\\ y & 0 \end{smallmatrix}\right)]$ or of type $[\left(\begin{smallmatrix}x & 0\\ y & x \end{smallmatrix}\right)].$ In the first case $[\left(\begin{smallmatrix}x & 0\\ y & ty  \end{smallmatrix}\right)]$ is family,
whose members are G-equivalent to points in $Z_R$ for $t\ne 0.$
In the second case the members of the family $[\left(\begin{smallmatrix}x & 
t^2y \\ y & x \end{smallmatrix}\right)]$ for $t\ne 0$ are G-equivalent to $[\left(\begin{smallmatrix}x+ty & 0\\ 0 & x-ty  \end{smallmatrix}\right)]$ belonging 
also to $Z_R.$
This proves $Y\subset\overline{GZ}_R$ and thus (1) and (2).
The lengthy but elementary proof of (3) and (4) by use of local 
coordinates for the Grassmannian and its blow-up is omitted here. 
\end{proof}

{\bf Remark:} The set $\overline{GZ}_R\smallsetminus GZ_R$ consists
entirely of the orbits of the unstable points  
$([\left(\begin{smallmatrix}x & 0\\
y & 0 \end{smallmatrix}\right)], [x\wedge y, 0, 0])$ and 
$([\left(\begin{smallmatrix}x & 0\\
y & x \end{smallmatrix}\right)], [\xi, 0, 0]).$

\begin{sub}\label{intersEG}{\bf Lemma:}
$E_G\smallsetminus E_G^{ss}\subset E_G\cap \overline{GZ}_R\subset 
E_G\cap\tilde{H}_1=E_G\smallsetminus E_G^{s}$\\
and these sets are of dimension 4,5,6 respectively.
\end{sub}

\begin{proof} When a point $p\in E_G$ is unstable, it is in the orbit
of a point $q=([\left(\begin{smallmatrix}x & 0\\
0 & x \end{smallmatrix}\right)], [\xi, 0, 0])$ and then $p\in\overline{GZ}_R.$ 
Such points have a 2-dimensional stabilizer $G_q$ and then 
$E_G\smallsetminus E_G^{ss}$ is parametrized by $P(Q)\times G/G_q,$
where $Q$ is the tautological quotient bundle on $P(V).$ Hence 
$E_G\smallsetminus E_G^{ss}$ is 4-dimensional.
The points in $E_G\cap \overline{GZ}_R$ are of type
$([\left(\begin{smallmatrix}x & 0\\
0 & x \end{smallmatrix}\right)], [a\xi, b\xi, c\xi,])$ with $\xi = x\wedge u$
and $u\in V/k.x$. Therefore there is a surjective morphism 
$P(Q)\times \P_2\to E_G\cap \overline{GZ}_R$ which is generically injective.
Hence $\dim (E_G\cap \overline{GZ}_R)=5.$ Finally $E_G\cap\tilde{H}_1$
is an intersection of hypersurfaces and so of dimension 6.
\end{proof}

The condition in Lemma \ref{gzr} for points in $\overline{GZ}_R$ is
equivalent to the vanishing of\\ 
$\xi\wedge\omega,\; \xi\wedge\eta,\; \omega\wedge\eta.$ Moreover, the 
homomorphism $(\xi,\; \omega,\; \eta)\mapsto (\xi\wedge\omega,\;\xi\wedge\eta,\; \omega\wedge\eta)$ describes the canonical wedge map
$$\Hom((S^2k^2)^*, \wedge^2V)\to \Hom(\wedge^2(S^2k^2)^*, \wedge^2\wedge^2V),$$
and this is $G$-equivariant, explicitly described by
\begin{equation}\label{action}(\xi,\; \omega,\; \eta)
\left(\begin{matrix} \alpha^2 & 2\alpha\beta & \beta^2\\
\alpha\gamma & \alpha\delta+\beta\gamma & \beta\delta\\
\gamma^2 & 2\gamma\delta & \delta^2  \end{matrix}\right)\mapsto
(\xi\wedge\omega,\;\xi\wedge\eta,\; \omega\wedge\eta)
\left(\begin{matrix} \alpha^2 & \alpha\beta & \beta^2\\
2\alpha\gamma & \alpha\delta+\beta\gamma & 2\beta\delta\\
\gamma^2 & \gamma\delta & \delta^2  \end{matrix}\right).
\end{equation}

So the map
$$\tilde{X}\smallsetminus\overline{GZ}_R\overset{\Phi}\lra P(\wedge^2(S^2k^2)\otimes\wedge^2\wedge^2V)\cong P(k^3\otimes V),$$
given by $p\to [\xi\wedge\omega,\; \xi\wedge\eta,\; \omega\wedge\eta]$
is well-defined and $G$-equivariant and the components of this map
generate the ideal sheaf of $\overline{GZ}_R$

{\bf The second Kirwan blow-up} can now be defined as the blow-up of
$\tilde{X}$ along $\overline{GZ}_R$:
$$Y:=Bl_{\overline{GZ}_R}(\tilde{X})\overset{\pi}\lra\tilde{X}$$
It is simultanously the {\bf closure of the graph of $\Phi.$}
By the smoothness of the ingridients, $Y$ is smooth. Moreover,
$Y\subset\tilde{X}\times P(k^3\otimes V)$ is acted on by $G$ and  
the projection $Y\to\tilde{X}$ is $G$-equivariant according to 
formula (3). We let $E_R$ denote the exeptional
divisor. 

{\bf Remark:} The condition for $\overline{GZ}_R$ says that the second
components of its points are of type $[a\xi, b\xi, c\xi]=(a,b,c)\otimes\xi$
in $P(S^2k^2\otimes\wedge^2V).$ This means that $\overline{GZ}_R$ is
the pull back of the Segre variety $S=P(S^2k^2)\times P(\wedge^2V)$
in $P(S^2k^2\otimes \wedge^2V).$ It follows that also the blow-up 
$Bl_{\overline{GZ}_R}(\tilde{X})$ is the pull back of the blow-up
of  $P(S^2k^2\otimes \wedge^2V)$ along the Segre variety $S.$

\begin{sub}\label{stab2}{\bf Stability in $Y$:}\rm\quad
By definition $Y$ is embedded in 
$\tilde{X}\times P(\wedge^2(S^2k^2)\otimes\wedge^2\wedge^2V).$
Combined with Segre embeddings we have
$$Y\subset P((\wedge^2k^2\otimes S^2V\oplus S^2k^2\otimes \wedge^2V)\otimes (S^2k^2\otimes \wedge^2V)\otimes (\wedge^2(S^2k^2)\otimes\wedge^2\wedge^2V)).$$
Then using the Mumford criterion and considering the weights of 
1-parameter subgroups one can derive: 

\begin{enumerate}
\item[(i)] Points in $Y$ over points in $GZ_R$ are stable.
\item[(ii)] Points in $Y$ over stable points in $\tilde{X}$ are stable.
\item[(iii)] Points in $Y$ over unstable points in $\tilde{X}$ are unstable.
\item[(iv)] Properly semistable points in $\tilde{X}^{ss}\smallsetminus GZ_R$
become unstable in $Y$.
\item[(v)] Every semistable point in $Y$ is stable.
\end{enumerate}
Remark: One can as well show that the stabilizer of any semistable point
in $Y$ is finite.
\end{sub}

The $G$-equivariant morphism $\pi$ induces a surjective $G$-equivariant
morphism $Y^s\lra\tilde{X}^{ss}$ and thus a surjective morphism of
the good quotients
$$\tau : Y^s/G\lra\tilde{X}^{ss}//G= \widetilde{P(S^2V)}$$
with surjective restriction
$$\tilde{\Sigma}_1:=E^s_R/G\lra GZ_R//G=\Sigma_1,$$
whereas 
$$Y^s/G\smallsetminus\tilde{\Sigma}_1\overset{\approx}\lra\tilde{X}^{ss}//G\smallsetminus\Sigma_1$$
must be an isomorphism because  $\pi$ is an isomorphism outside $E_R.$ 
Moreover, because $Y^s\overset{\pi}\lra\tilde{X}^{ss}$ is a blow up,
also the induced morphism $\tau$ is a blow-up along the divisor ${\Sigma}_1.$
Hence the 

\begin{sub}\label{fibres}{\bf Proposition:}
$\tau : Y^s/G\lra\tilde{X}^{ss}//G = \widetilde{P(S^2V)}$
is an isomorphism.
\end{sub}

Remark: While the second Kirwan blow-up has no effect on the quotient,
it describes $\widetilde{P(S^2V)}$ as a geometric quotient, so that
non-isomorphic S-equivalent limit sheaves w.r.t. the parameter space
$Y^s$ are excluded. This is needed for the construction of 
families which include admissible tree bundles because 
S-equivalence for tree bundles is not defined.
\end{section}

\section{Families including tree bundles}

In this section the construction of families of sheaves,
including all admissible tree bundles for the tree compactification of 
$M^b(0,2)$, 
will be sketched in two steps. In step one we construct such a family
over the base space $\tilde{X}^{ss}$.

Firstly we recall the presentation of the semi-universal family for the 
Gieseker-Maruyama space $M(2;0,2)$. Let
$0\to\ku\to k^2\otimes\ko_X\to\kq\to 0$ be the tautological sequence on the 
Grassmannian $X=G_2(k^2\otimes V).$ As in formula (\ref{beilin}) there are two
such equivalent presentations. The second is the exact sequence over 
$X^{ss}\times P$
\begin{equation}\label{univ}
0\to k^2\otimes\ko_X\boxtimes\ko_P(-2)\lra\kq\boxtimes\ko_P(-1)\to\kf\to 0.
\end{equation}

Recall from \ref{remstab} that $H_1^{ss}\subset X^{ss}$ is the hypersurface
of points $[A]$ for which $\det(A)$ decomposes, i.e. the inverse image
of $\Sigma_1$, and that $H_0^{ss}\subset H_1^{ss} $
is the subvariety where $\det(A)$ is a square. Let now 
$S_1\subset X^{ss}\times P$ be the subvariety of points $([A],[v])$
for which $v$ divides $\det(A)$, and $S_0\subset S_1$ where
$\det(A)=v^2.$ Then $S_1$ is 7-dimensional and 2:1 over 
$H_1^{ss}\smallsetminus H_0^{ss}.$

It follows that $\kf$ is locally free on $X^{ss}\times P\smallsetminus S_1$
whose restriction to fibres over $X^{ss}\smallsetminus H_1$
are the vector bundles in $M^b(0,2),$ whereas the sheaves over points 
in $H_1$ become the semistable sheaves in the boundary of $M^b(0,2).$

{\bf Notice} however that the sheaf $\kf$ restricted to $\{p\}\times P$ 
may be singular only in one of the points of $S_1$
over $p$, see the Notice before \ref{typ1}.

\begin{sub}\label{step1}{\bf First step:}\rm

Let now $\tilde{X}^{ss}\times P\xra{\alpha} X^{ss}\times P$ be the map
$\phi=\alpha\times\id$, where $\alpha$ is the blow-up map 
of section 4, and consider the lifted family $\F=\phi^*\kf.$
Then $\F$ is locally free over the inverse image of 
$X^{ss}\smallsetminus E_G^{ss}\cup\tilde{H}_1^{ss}$

Analogously to $S_0$ and $S_1$, let then $\tilde{S}_0$ the set of points 
$(p,[v])\in\tilde{X}^{ss}\times P$ over $E_G^{ss}$  where
$\det(A)=v^2$, and let similarly  
$\tilde{S}_1\subset\tilde{X}^{ss}\times P$ be the set of points over 
$\tilde{H}_1^{ss}$ where $v$ is a factor of $\det(A)$.
Then $\F$ is locally free outside $\tilde{S}_0\cup\tilde{S}_1,$
$\tilde{S}_0$ is mapped 1:1 to $E_G^{ss}$ and the map 
$\tilde{S}_1\smallsetminus\tilde{S}_0\to\tilde{H}_1^{ss}\smallsetminus E_G^{ss}$ is 2:1. 

Consider now the blow-up $Z\xra{\sigma_0}\tilde{X}^{ss}\times P$ along 
$\tilde{S}_0$
and let $D$ denote the exceptional divisor. Let 
$$\bkf:=\sigma_0^*\F/torsion$$
be the torsion free pullback on $Z.$
Now the situation of the families $\bkf$ and $\F$ restricted to the 
open subset 
$\tilde{X}^{ss}\smallsetminus\tilde{H}_1^{ss}$ of the base  
is the higher dimensional analog to that  
of the families over the curve $C$ in section \ref{typ1},
with $0\in C$ replaced by the divisor 
$E_G^{ss}\subset\tilde{X}^{ss}\smallsetminus\tilde{H}_1^{ss}.$

Moreover, one can compare the two situations
by considering a curve $C\subset\tilde{X}^{ss}$ transversal to $E_G^{ss}$ in a 
point $p\not\in\tilde{H}_1^{ss}.$ Then the blow-up $Z_C$ of $C\times P$ in the 
point $(p,q)\in \tilde{S}_0$ can be 
identified with the restriction of $Z$ to $C$. Moreover, by flatness, the 
sheaves $\bkf_C$ and $\bF_C$ on $Z_C$ from \ref{typ1} can be identified
with the restrictions of $\bkf$ and $\bF$ to $Z|C.$ 
Because $\bkf_C$ is locally free on $Z_C$, it follows that $\bkf$ is locally
free in a neighborhood of the fibre $Z_p$ of $Z$ over $p$.
Finally, because the fibre $Z_p$ is the union of the blow-up of $P$ at $q$
and the restriction $D_p$ of exceptional divisor $D$, the sheaf  $\bkf|Z_p$
is a tree bundle on $Z_p$. In order to obtain the correct Chern classes,
we have to replace $\bkf$ by its twist $\bkf(D)$ as in \ref{typ1}, which
is also compatible with the restriction. It has been shown:

\begin{subsub}\label{trees1}{\bf Proposition:}
With the notation above, the family $\bkf$ is a family of tree bundles
over the restricted base variety 
$\tilde{X}^{ss}\smallsetminus\tilde{H}_1^{ss}.$\\
If $p\in\tilde{X}^{ss}\smallsetminus\tilde{H}_1^{ss}\cup E_G^{ss},$ then 
$\bkf|Z_p,\; Z_p=P,$ is a bundle in $M^b(0,2)$.\\
If 
$p\in E_G^{ss}\smallsetminus\tilde{H}_1^{ss},$
then  $\bkf|Z_p,$ where $Z_p=\tilde{P}\cup D_p,\; D_p\cong \P_2,$ 
is a tree bundle with
$\bkf|\tilde{P}\cong 2\ko_{\tilde{P}}$ and $\bkf|D_p\in M^b_{D_p}(0,2).$
\end{subsub}

For the fibres over points in $\tilde{H}_1^{ss}$ we have:
\begin{subsub}\label{trees2}{\bf Lemma:}
Let $\hat{S}_1$ be the proper transform of $\tilde{S}_1$ in $Z.$
Then $\hat{S}_1\to \tilde{H}_1^{ss}$ is 2:1.
\end{subsub}
Remark: For a point $p\in E_G^{ss}\cap\tilde{H}_1^{ss}$ the two points
of $\hat{S}_1$ over $p$ will be contained in the fibre $D_p\cong\P_2$
of the exceptional divisor $D.$ By the previous,
$\bkf$ is locally free on $Z\smallsetminus\hat{S}_1.$

\begin{proof} 
The method of proof is again by restriction to transversal curves:
Let $q\in\tilde{S}_0\cap\tilde{S}_1\subset\tilde{X}^{ss}\times P$ and
$p\in E_G^{ss}\cap\tilde{H}_1^{ss}$ its image. Then $p$ has the components
$[\begin{smallmatrix} x & 0\\ 0 & x \end{smallmatrix}]$ 
and $[a\xi,b\xi,c\xi]$ with $b^2-4ac\ne 0$ and $\xi=x\wedge y$ for some
$y\in V$. Then
$$
p(t):=([\begin{smallmatrix} x & 0\\ 0 & x+ty \end{smallmatrix}], 
[a\xi,b\xi,c\xi])
$$
is a 1-parameter family in $GZ_R\subset\tilde{H}_1^{ss}$ defining a normal 
direction
to $E_G^{ss}$ at $p$. Let $C$ denote the image of $p(t)$ for small $t.$
For $t\ne 0$ the points [x] and [x+ty] define then sections of 
$\tilde{S}_1|C\smallsetminus\{0\},$ which fill this subset. Because $q$
is the only point in $\tilde{S}_0$ over $p$, $q\in\tilde{S}_1|C$, the closure
of $\tilde{S}_1|C\smallsetminus\{0\}.$
Let now $S_C:=\tilde{S}_1|C\subset C\times P$ and consider the blow-ups
$$ Bl_q(C\times P)\subset Bl_{\tilde{S}_0}(\tilde{X}^{ss})=Z$$
as the proper transform. Then the restriction $\hat{S}_1|C$ of the
proper transform $\hat{S}_1$ can be identified with the proper transform
of $S_C$ in $Bl_q(C\times P).$ This situation corresponds to the figure in
\ref{typ1} with the two sections [x] and [x+ty] added. Then the proper
transforms of these linear sections do not meet on the exceptional divisor
$D_p.$ Hence also $\hat{S}_1\cap D_p$ consists of two different points.
\end{proof}
\end{sub}

\begin{sub}\label{step2}{\bf Second step:}\rm

By the above, $\bkf$ is locally free on $Z\smallsetminus\hat{S}_1$
and one could try to construct the tree bundles over $\tilde{H}_1^{ss}$ 
by directly blowing up $Z$ along $\hat{S}_1$ and modifying the lifted
sheaf. However, over points $p\in\tilde{H}_1^{ss}\smallsetminus GZ_R$
the sheaf $\bkf|Z_p$ has only one singular point and is not stable,
see the remark at the beginning of this section 6. Secondly, 
$\tilde{H}_1^{ss}\smallsetminus GZ_R$ consists only of non-closed orbits.
On the other hand the orbits in $GZ_R$ are closed and for $p\in GZ_R$
the two points of $\hat{S}_1$ are the singular points of  $\bkf|Z_p.$

Now this insufficiency can be eliminated by using the second Kirwan
blow-up $Y\to\tilde{X}$ and pulling the pair $(Z,\bkf)$ back to
$Y^s$. After this the points of $\tilde{H}_1^{ss}\smallsetminus GZ_R$
become unstable and can be neglected, and $\hat{S}_1| GZ_R$ is the reasonable
locus to be blown up. Therefore, let
\[
\xymatrix{Z_Y\ar[d]\ar[r] & Z\ar[d]\\
Y^s\ar[r] & \tilde{X}^{ss},
}
\]
be the pull back of $Z$ and let $\bkf_Y$ be the lift of $\bkf$ to $Z_Y.$
The situation of the pair $(Z_Y,\bkf_Y)$ is now the relative version
of the situation in \ref{typ2} before using a double cover.

\begin{subsub}\label{onY}{\bf Properties of $(Z_Y,\bkf_Y)$:}\rm\\
Let $E^s_R$ denote the exceptional divisor of
$Y^s$ over $GZ_R,$ see section 5, let $E_{G,Y}$ denote the proper transform 
of $E_G,$ and let $D_Y$ be the pull back of $D$ in $Z$. 
Then $\bkf_Y$ is singular exactly along the pull back
$S_{1,Y}$ of $\hat{S}_1$ and $S_{1,Y}$ is 2:1 over $E^s_R$ everywhere
by \ref{trees2}. For points $p$ in $E^s_R\smallsetminus E_{G,Y}$,  
the two points of $S_{1,Y}$ over $p$ will be in the fibre
$Z_{Y,p}\cong P,$ but for points $p$ in $E^s_R\cap E_{G,Y}$,
the two points of $S_{1,Y}$ over $p$ will be in the fibre $D_{Y,p}$ of 
$D_Y$. 
\end{subsub}

{\bf Remark:} The variety $Z_Y$ may also be obtained as the blow-up
of the variety $S_{0,Y}\subset Y^s\times P$ over $E_{G,Y},$ defined
as $\tilde{S_0}$ over $E_G.$

In order to construct a family of tree bundles in this new relative 
situation, $Z_Y$ has to be blown
up along $S_{1,Y}$ as in the case \ref{typ2}. Then the torsionfree
pull back of $\bkf_Y$ would give a family of tree bundles parametrized
along $E^s_R$. But as in \ref{typ2} these tree bundles would not be
admissible as defined in \ref{trvb}. In analogy to \ref{typ2}
one would have to use a double cover of $Y^s$ which is branched exatly over 
$E^s_R$ in order to construct admissible tree bundles. However, such a
double cover may not exist globally. But one could consider such
local covers $U\to Y^s$ over affine open parts. Then we have 
Cartesian diagrams
\[
\xymatrix{W_U\ar[r]^\tau & Z_U\ar[d]\ar[r]^{g} & Z_Y\ar[d]\\
  & U\times P\ar[d]\ar[r] & Y^s\times P\ar[d]\\
  & U \ar[r]^f & Y^s,
}
\]
where $\tau$ is the blow-up of $Z_U$ along the subvariety $S_U=g^*S_{1,Y}.$
This is the subvariety where $\bkf_U:=g^*\bkf_Y$ is not locally free.
By the previous, it is 2:1 over the branch locus $B:=f^*E^s_R\subset U.$
Consider then the sheaf
$$\bke:=\tau^*\bkf_U/torsion.$$
One can show as in the curve case that $\bke$ is flat over $U.$

Now one can argue as in \ref{step1} using curves $C$ which are transversal
to $B:$ There is an elementary transform $\bke^\prime$ of $\bke$ on $W_U$
with transformation support over $B$ which is locally free on 
$W_U$. Then 
$\bke^\prime$ is a family of tree bundles, whose fibres over points
in $U\smallsetminus B$ are the same as for points in
$Y^s\smallsetminus E^s_R$ or in $\tilde{X}^{ss}\smallsetminus GZ_R.$
After twisting with the exceptional divisor in $W_U$, we may finally 
assume that $\bke^\prime$ is a family of admissible tree bundles with prescribed
Chern classes. Hence the  

\begin{subsub}\label{final}{\bf Proposition:} For any 2:1 cover 
$U\xra{f}Y^s$ of an affine open subset of $Y^s, $branched exactly along 
$E^s_R$, the following holds: 
\begin{enumerate}
\item[(i)] For points $p$ in $U\smallsetminus f^*\tilde{E}_G\cup B$
the bundle 
$\bke^{\prime}_p$ is a member of $M^b_{W_{U,p}}(0,2),$ where $W_{U,p}\cong\P_2$.
\item[(ii)] For points $p$ in $f^*\tilde{E}_G\smallsetminus B$
the bundle $\bke^{\prime}_p$ is of the type described in \ref{typ1}. 
\item[(iii)] For points $p$ in $B\smallsetminus f^*\tilde{E}_G$
the bundle $\bke^{\prime}_p$ is of the type described in \ref{typ2}.
\item[(iv)] For points $p$ in $B\cap f^*\tilde{E}_G$
the bundle $\bke^{\prime}_p$ is of the type described in \ref{typ3}.
\end{enumerate}
\end{subsub}

The families of tree bundles so constructed may not descend to a global
family over the Kirwan blow-up $\tilde{M}_2\cong\widetilde{P(S^2V)}$ 
of $M(2;0,2)\cong P(S^2V)$ because the automorphism groups of the
tree bundles include automorphisms of the supporting surfaces, see
\ref{trvb}. However, delicately, their isomorphism classes  
are determined precisely by the points of $\tilde{M}_2$:

\begin{subsub}\label{final2}{\bf Proposition:}\quad The set of points
of $\tilde{M}_2$
is the set isomorphism classes of the tree bundles constructed above. 
In particular, let as above 
$\Sigma_2\subset\widetilde{P(S^2V)}$ be the blow-up of the Veronese
surface in $P(S^2V)$ and $\Sigma_1\subset\widetilde{P(S^2V)}$ the proper 
transform of the subvariety of decomposable conics. Then

\begin{enumerate}
\item[(i)] $\tilde{M}_2\smallsetminus\Sigma_1\cup\Sigma_2=M^b(0,2)$ is the set
the isomorphism classes of the (stable) bundles in $M(2;0,2).$
\item[(ii)] The set $\Sigma_2\smallsetminus\Sigma_1$ is
the set of isomorphism classes of limit tree bundles of type 1 
described in \ref{typ1}.
\item[iii)] The set $\Sigma_1\smallsetminus\Sigma_2$ is the set of
isomorphism classes of limit tree bundles of type 2 described in 
\ref{typ2}.
\item[(iv)] The set $\Sigma_1\cap\Sigma_2$ is the set of
isomorphism classes of limit tree bundles of type 3 described in 
\ref{typ3}.
\end{enumerate}
\end{subsub}

\begin{proof} There is nothing to proof for (i). For the proof of (ii),
recall that $\Sigma_2\smallsetminus\Sigma_1$ is the geometric quotient
of the open part $E_G^s\subset\tilde{X}^{s}$ of the exceptional divisor
$E_G$ whose points are of type
$$p=([\left(\begin{smallmatrix} x & 0 \\
 0 & x \end{smallmatrix}\right)], [\xi,\; \omega,\; \eta]),$$
where $\omega^2-4\xi\eta$ decomposes into two different factors and
$\xi,\; \omega,\; \eta\in x\wedge V.$ By \ref{remstab} we may assume that
$\omega =0.$
The two factors $\xi,\;\eta$ determine two lines in P=P(V) through $[x]$,
see \ref{dual}. Now the fibre $Z_p$ is a union $\tilde{P}(x)\cup D_p$,
where $\tilde{P}(x)$ is the blow-up of $P$ at $[x]$ and $D_p\cong\P_2.$
Then the two lines in $P$ determine two points $q_1,\; q_2$ on the 
exceptional line $\ell_p=\tilde{P}(x)\cap D_p.$
Let now $\bkf$ on $Z$ be the sheaf constructed in \ref{step1}. By
\ref{trees1} $\bkf|Z_p$ has the restrictions 
$\bkf|\tilde{P}(x)\cong 2\ko_{\tilde{P}(x)}$ and $\bkf|D_p\in M^b_{D_p}(0,2).$
So $\bkf|D$ corresponds to its smooth conic of jumping lines in the dual
plane $D_p^*$ or to the dual conic $\Gamma_p\subset D_p$ of the latter.

{\bf Claim:} The conic  $\Gamma_p$ meets the line $\ell_p$ in the two points 
$q_1,\; q_2.$

In addition, there is the following elementary 

\begin{subsub}\label{lemii}{\bf Lemma:} Let $\ell$ be the line through two 
points $a_1,\; a_2\in\P_2$
and let $\Aut_\ell(\P_2)$ be the subgroup of the group of automorphisms
of $\P_2$ which fixes the points of $\ell$. Then $Aut_\ell(\P_2)$
acts transitively on the set of non-degenerate conics through $a_1,\; a_2$.
\end{subsub}

If the claim is verified, the Lemma implies that the isomorphism class
of  $\bkf|D_p$ and then also of $\bkf|Z_p$ only depends on the two points
$q_1,\; q_2$, which are determined by the point $p.$ Then the 
isomorphism class of $\bkf|Z_p$ also depends only on the image $[p]$
of $p$ in the quotient $\tilde{M}_2$, which proves (ii).

In order to prove the claim we use again 1-parameter degenerations with
limit point $p$ which are transversal to $E_G^s$. 

For that we may assume that 
$$p=([\left(\begin{smallmatrix} e_0 & 0 \\
 0 & e_0 \end{smallmatrix}\right)], [e_0\wedge e_1,\; 0,\; e_0\wedge e_2]),$$ 
where $e_0,e_1,e_2$ form a basis of $V$, and that the first component
of the 1-parameter family is given by
$$A(t)=\left(\begin{smallmatrix} e_0 & 0 \\
 0 & e_0 \end{smallmatrix}\right) + t\left(\begin{smallmatrix}x & x^\prime\\
y & y^\prime   \end{smallmatrix}\right).$$
with $t$ in a neighborhood $C$ of $0\in \A^1(k).$ This is a smooth curve in
$X^{ss}$ and its lift to $\tilde{X}^{ss}$ is transversal to $E_G^s$ and has 
second component
$$[e_0\wedge y+t\xi,\; e_0\wedge (y^\prime-x)+t\omega,\; -e_0\wedge x^\prime+t\eta],$$
where $(\xi,\; \omega,\; \eta)=(x\wedge y,\; x\wedge y^\prime+x^\prime\wedge y,\; x^\prime\wedge y^\prime).$ 
Because $p$ is supposed to be the limit at $t=0$,
we may assume, up to a scalar factor, that the components of the vectors satisfy
$$y_1=1,\; y_2=0,\; y_1^\prime=x_1,\; y_2^\prime=x_2,\; x_1^\prime=0,\; x_2^\prime= -1$$  
In addition we replace the basis 
$e_1\wedge e_2,\; -e_0\wedge e_2,\; e_0\wedge e_1$  of $\wedge^2V$ by 
the basis $z_0,z_1,z_2$ of $V^*$, dual to the basis $e_0,e_1,e_2$ of $V$.
Then the second component of $p(t)$ reads
$$[z_2-tx_2z_0+t\xi^\prime,tz_0+t\omega^\prime,-z_1+x_1z_0+t\eta^\prime],$$
where $\xi^\prime,\; \omega^\prime,\; \eta^\prime\in \Span(z_1,z_2)$.

Let now $Z_C$ be the restriction  of $Z$ to $C$. Then $Z_C$ can considered
the blow up of $C\times P$ at $(0,[e_0])$ as a proper transform and 
$\bkf|D_p$ can be computed as in \ref{typ1}, as well as its conic 
$\Gamma_p\subset D_p$. As $\bkf|Z_C$ is the torsion free pull back
of the sheaf on $C\times P$ defined by $A(t)$, its family of conics becomes 
the proper transform of the family 
$$Q=\{(tz_0+t\omega^\prime)^2-4(z_2-tx_2z_0+t\xi^\prime)(-z_1+tx_1z_0+t\eta^\prime)=0\},$$
whose fibres for $t\ne 0$ are the conics of $\bkf|\{t\}\times P$, 
c.f. \ref{dual}.
This proper transform is obtained by substituting the forms 
$tz_0, z_1, z_2$ by $u_0,u_1,u_2$, which are the coordinate forms of 
$D_p\cong\P_2$, see \ref{typ1}. So the proper transform $\tilde{Q}$ of $Q$ 
is defined by the equation
$$(u_0+t\omega^\prime)^2-4(u_2-x_2u_0+t\xi^\prime)(-z_1+x_1u_0+t\eta^\prime),$$
where now $\xi^\prime,\; \omega^\prime,\; \eta^\prime\in \Span(u_1,u_2)$.
For $t=0$ the conic $\Gamma_p$ of $\bkf|D_p$ has the equation 
$u_0^2-4(u_2-x_2u_0)(-u_1+x_1u_0)$. Because the line $\ell_p$ and 
$\tilde{P([e_0])}$ is given by $u_0=0$, $\Gamma_p$ meets $\ell_p$
in the two points $q_1, q_2$ with equation $u_1u_2.$ 
This proves the claim and thus (ii) of proposition \ref{final2}.

For the proof of (iii) let a point in $\Sigma_1\smallsetminus\Sigma_2$
be the image of a point $p\in E_R^s\subset Y^s.$ We may assume that
its image $\bar{p}\in GZ_R$ under the second blow up has the components
$$\left(\begin{smallmatrix} x_1 & 0\\ 0 & x_2  \end{smallmatrix}\right), [0,\; x_1\wedge x_2,\; 0].$$
Under an auxiliary blow-up $W_U$ as in \ref{onY},
the fibre $W_{U,p}$ is isomorphic to $\tilde{P}(x_1,x_2)\cup D_1\cup D_2,$
where $\tilde{P}(x_1,x_2)$ is the blow-up of $P$ at $x_1,x_2$, $D_i\cong\P_2$,
containing the exceptional lines $\ell_i$ of $\tilde{P}(x_1,x_2).$ 
Moreover, $W_{U,p}$ is determined by the data of the point $p$ or its 
image in $\Sigma_1\smallsetminus\Sigma_2$ up to isomorphism.
By \ref{final}, (iii), the tree bundle $\bke^{\prime}_p=\bke^{\prime}|W_{U,p}$
is trivial on $\tilde{P}(x_1,x_2)$ and restricts to bundles $\bke^{\prime}_i$
on $D_i$ with Chern classes $c_1=0, c_2=2.$ By the following 
Lemma \ref{lemiii}
the isomorphism class of each $\bke^{\prime}_i$ corresponds uniquely to a point 
$q_i\in D_i\smallsetminus \ell_i$. Since the group $\Aut_{\ell_i}(D_i)$
acts transitively on $D_i\smallsetminus \ell_i$, see \ref{lemii},
these isomorphism classes are uniquely determined by $W_{U,p}$,
and finally determined by the point $[p]\in\Sigma_1\smallsetminus\Sigma_2$,
because the automorphisms of $W_{U,p}$ must be identities on 
$\tilde{P}(x_1,x_2).$ 

\begin{subsub}\label{lemiii}{\bf Lemma:}\quad Let $\ell\subset\P_2=P$ be a line.
Then the moduli space $M_\ell(0,1)$ of isomorphism classes of rank-2 vector 
bundles on $\P_2$ which are trivial on $\ell$ with Chern classes $c_1=0, c_2=1$
can be identified with the set $\P_2\smallsetminus\ell.$
\end{subsub}

Proof of the Lemma: Let $\ell$ have the equation $z_0$ and let 
$a=[a_0,a_1,a_2]\in\P_2\smallsetminus\ell.$
Let 
$$B=\left(\begin{matrix} z_0 & z_1 & z_2 & 0\\
a_0 & a_1 & a_2 & z_0 \end{matrix}\right),$$ and define $\ke(a)$ as cokernel
in the sequence
$$0\to\ko_P(-2)\oplus\ko_P(-1)\xra{B}3\ko_P(-1)\oplus\ko_P\to\ke(a)\to 0.$$
Then the class of $\ke(a)$ belongs to $M_\ell(0,1).$ Conversely, given
any $\ke$ in $M_\ell(0,1)$, it is well known that $\ke$ is an elementary
transform of the twisted tangent bundle $\kt_P(-2)$ with exact extension
sequence
$$0\to\kt_P(-2)\to\ke\to\ko_\ell\to 0.$$
From that we get a resolution matrix $B$ of $\ke$ as above. In that,
$(a_0,a_1,a_2)$ represents the extension class and $[a_0,a_1,a_2]$
the isomorphism class of $\ke$.

This completes the proof of (iii). The proof of (iv) is analogous to
that of (iii). In this case a point $p\in B\cap f^*\tilde{E}_G$ 
or $p\in E^s_R\cap\tilde{E}_G$ over a point in 
in $\Sigma_1\cap\Sigma_2$ can be supposed to have as components 
$$[\left(\begin{smallmatrix} x & 0\\ 0 & x  \end{smallmatrix}\right)], 
[a\xi,\; b\xi,\; c\xi],[u,w,v]$$
with $\xi\in x\wedge V$ and $b^2-ac\ne 0.$
Then $W_{U,p}$ as a fibre of the blow-up is isomorphic to 
$\tilde{P}([x])\cup\tilde{D}_0(p_1,p_2)\cup D_1\cup D_2,$
where $\tilde{D}_0(p_1,p_2)$ is the blow-up at two points of a plane $D_0$
which contains the exceptional line $\ell_0$ of $\tilde{P}([x])$,
and where $D_i$ are again planes containing the two exceptional lines
$\ell_i$ of $\tilde{D}_0(p_1,p_2).$ Then $W_{U,p}$ depends only on the geometry
and the point $p$ up to isomorphism.
Now the tree bundle $\bke^{\prime}$ on $W_U$ of \ref{final} is trivial
on $\tilde{P}([x])$ and $\tilde{D}_0(p_1,p_2)$, whereas $\bke^{\prime}|D_i$
has Chern classes $c_1=0, c_2=1.$ It follows again from \ref{lemiii}
that the isomorphism classes of $\bke^{\prime}|D_i$ are unique, and then
that $\bke^{\prime}|W_{U,p}$ is uniquely determined
because the automorphisms of $W_{U,p}$ must be identities on the 
components $\tilde{P}([x])$ and $\tilde{D}_0(p_1,p_2)$. This proves (iv)
of the proposition.
\end{proof}

%{\bf Remark:} In order to obtain a global 2:1 cover of $Y^s$ which
%is branched exactly along $E^s_R$, one could replace $Y^s$ by the blow-up
%of $\tilde{X}^{ss}$ with respect to the square $\ki^2_R$ of the 
%ideal sheaf $\ki_R$ of $GZ_R$, which is isomorphic to $Y^s$ so that
%the invertible sheaf of its exceptional divisor is a square.
\end{sub}

\begin{sub}{\bf The stack:}\rm

By the above construction of families of tree bundles a global
family of such bundles could not be obtained. Instead, we have
families of tree bundles on local 2:1 covers of the parameter space
$Y^s$. These are forming an obvious moduli stack over the category
of such open covers. It is plausible to claim that this is a 
Deligne-Mumford stack which is corepresented by 
$Y^s/G\cong\widetilde{P(S^2V)}.$ 

In \cite{MTT} global families of limit tree bundles 
of stable rank-2 vector bundles on surfaces have been constructed by 
other abstract procedures, which led to algebraic spaces as moduli spaces.
The question of their relation to the above stack 
being open at present.
\end{sub}
%\newpage

\end{document}